\documentclass[11pt]{article}

\usepackage{amsmath}
\usepackage{amssymb}
\usepackage{amscd}
\usepackage[latin1]{inputenc}
\usepackage{graphics}

\newcommand{\Keler}             {K\"{a}hler }

\newcommand{\OO}{\mathcal{O}}
\newcommand{\meno}{^{-1}}
\newcommand{\Zeta}{{\mathbb{Z}}}
\newcommand{\Spin}{\operatorname{Spin}}
\newcommand{\Lam}{\Lambda}

\newcommand{\W}{\mathbb{W}}

\newcommand{\V} {U}

\newcommand{\spam}{\,\operatorname{span}}

\newcommand{\fund}{\varpi}

\newcommand{\pot}[1][]{{\, \Lambda^r#1}}

\newcommand{\chern}{\operatorname{c}}
\newcommand{\comp}{\circ}

\newcommand{\alfa}{\alpha}

\newcommand{\gh}{{K}}

\newcommand{\lieh}{\mathfrak{h}}
\newcommand{\liel}{\mathfrak{l}}
\newcommand{\lieg}{\mathfrak{g}}
\newcommand{\lieb}{\mathfrak{b}}
\newcommand{\liep}{\mathfrak{p}}
\newcommand{\liez}{\mathfrak{z}}
\newcommand{\lies}{\mathfrak{s}}

\newcommand{\m}{{\mathcal{K}(L)}} 
\newcommand{\la}{\lambda}

\newcommand{\enf}{\emph}

\newtheorem{teo}{Theorem}

\newtheorem{prop}[teo]{Proposition}
\newtheorem{cor}[teo]{Corollary}
\newtheorem{lemma}[teo]{Lemma}

\newcommand{\Inn}{\operatorname {Inn}}

\newcommand{\restr}[1]          {\phantom{}_{\text{\raisebox{.4ex}{$|$}}#1}}

\newcommand{\w}{\mathbb{W}}

\newcommand{\Ad}{\operatorname{Ad}}
\newcommand{\ad}{{\operatorname{ad}}}

\newcommand{\Aut}{\operatorname{Aut}}
\newcommand{\Hom}{\operatorname{Hom}}
\newcommand{\SU} {\operatorname{SU}}
\newcommand{\su} {\mathfrak{su}}
\newcommand{\Sl}{\operatorname{SL}}
\newcommand{\Sp}{\operatorname{Sp}}

\newcommand{\dimo}[1][]         {\noindent\textbf{Proof#1}. }

\newcommand{\fine}{\hfill$\Box$\medskip \\ }

\newcommand{\ts}{\widetilde{S}}

\newcommand{\cinf}{C^\infty}

\newcommand{\Gl}{\operatorname{GL}}
\newcommand{\ra}{\rightarrow}
\newcommand{\C}{\mathbb{C}}
\newcommand{\KE}{{K\"ahler-Einstein} }
\newcommand{\om}{\omega}
\newcommand{\eps}{\varepsilon}
\renewcommand{\phi}{\varphi}

\newcommand{\scalar}{\langle \cdot \, , \cdot \rangle}

\newcommand{\basis}[1]{{{\boldsymbol{#1}}}}

\begin{document}

\title{
Homogeneous bundles and the first eigenvalue of symmetric
  spaces}

\author{Leonardo Biliotti, Alessandro Ghigi} \date{} \maketitle

{\abstract{
    \noindent In this note we prove the stability of the Gieseker
    point of an irreducible homogeneous bundle over a rational
    homogeneous space. As an application we get a sharp upper estimate
    for the first eigenvalue of the Laplacian of an arbitrary K\"ahler
    metric on a compact Hermitian symmetric spaces of ABCD--type.  }
  
}

\section{Introduction}

Let $X$ be a compact complex manifold and let $E$ be a holomorphic
vector bundle of rank $r$ over $X$.  The \enf{Gieseker point} of $E$
is the map
\begin{equation}
  \begin{gathered}
    T_E: \Lambda^r H^0(X,E) \longrightarrow H^0(X, \det E) \\
     \end{gathered}
\end{equation}
that sends an element $s_1\wedge \dots \wedge s_r \in \Lambda^r
H^0(X,E)$ to the section $ x \mapsto s_1(x) \wedge \cdots \wedge
s_r(x)$ of $H^0(X,\det E)$.  This map was first considered by Gieseker
in his work \cite{gieseker-vector-surfaces} in order to construct the
moduli space of vector bundles on a projective manifold.  He proved
that for the set of Gieseker stable bundles $E$ with fixed rank and
Chern classes on a polarised $(X,H)$ there is a uniform $k_0 $ such
that for $k>k_0$, $T_{E(k)} = T_{E \otimes H^{\otimes k}}$ is a stable
vector (in the sense of geometric invariant theory) with respect to
the action of $\Sl\bigl ( H^0(X,E(k)) \bigr )$ on $\Hom \bigl ( \Lambda^r
H^0(X,E(k)), H^0(X,\det E(k)) \bigr)$.

In this paper we consider the Gieseker point of homogeneous bundles
over rational homogeneous spaces.  Such bundles are known to be
Mumford-Takemoto stable \cite{ramanan-homogeneous},
\cite{umemura-homogeneous} and this implies they are Gieseker stable
\cite[p.191]{kobayashi-vector}. So we already know that after
twisting with a sufficiently ample line bundle their Gieseker point is
stable.  The interest here is in the stability of $T_E$ itself,
without allowing any twist. Our result is the following.
\begin{teo} \label{gieseker-intro} Let $E\ra X$ be an irreducible
  homogeneous vector bundle over a rational homogeneous space $X=G/P$.
  If $H^0(E) \neq 0$, then $T_E$ is stable.
\end{teo}
We give two proofs of this result. The first is algebraic and uses a
criterion of Luna for an orbit to be closed. This proof works over any
algebraically closed field of characteristic zero. The second proof
uses invariant metrics and relies on a result by Xiaowei Wang
\cite{wang-xiaowei-balance}.

Our interest for this result is connected with a problem in \Keler
geometry.  Consider a compact \Keler manifold $X$ and fix a \Keler
class $a \in H^2(X)$.  Bourguignon, Li and Yau
\cite{bourguignon-li-yau}, gave an upper bound for the first
eigenvalue of the Laplacian $\Delta_g : C^\infty(X) \ra \cinf(X)$
relative to any \Keler metric $g$ whose \Keler form $\om_g$ lies in
the class $a$. The bound depends on the numerical invariants ($h^0$
and degree) of a globally generated line bundle $L$ over $ X$.  To get
the best estimate one has to choose appropriately the bundle. As is
shown in \cite{bourguignon-li-yau}, if $X=\mathbb{P}^n$ and
$L=\OO_{\mathbb{P}^n}(1)$ one gets the upper bound $2$, which is
optimal since it is achieved by the Fubini-Study metric.

In paper \cite{arezzo-ghigi-loi} Arezzo, Loi and the second author,
generalised this result, by substituting a vector bundle $E$ to the
line bundle $L$. In this case one gets the same kind of estimate, but
the vector bundle $E$ must satisfy an additional condition, namely its
Gieseker $T_E$ point must be stable.  By this method one gets an upper
bound for $\la_1$ on the complex Grassmannian \cite{arezzo-ghigi-loi}.
Such a bound is optimal since it is achieved by the symmetric metric.
It is important to notice that if one twists the bundle $E$ by a
positive line bundle $H$, the estimate gotten from the twisted bundle
$E(k)$ is very rough. In fact the estimate blows up as $k \ra \infty$
(see \eqref{eq:def-J} below). So it is important to obtain some
information on the stability of $T_E$ without twisting the bundle $E$.

The main motivation for the present work was to extend the estimate
for $\la_1$ to other Hermitian symmetric spaces of the compact type
using appropriate homogeneous bundles.  We are able to prove the
following.
\begin{teo}\label{classical-intro}
  Let $X$ be a compact irreducible Hermitian symmetric space of
  ABCD-type.  Then
\begin{equation*}
  \la_1(X, g) \leq 2
\end{equation*}
for any \Keler metric $g$ whose \Keler class $\om_g$ lies in
$2\pi\chern_1(X)$. This bound is attained by the symmetric metric.
\end{teo}
It should be mentioned here that El Soufi and Ilias \cite[Rmk. 1,
p.96]{el-soufi-ilias-Pacific} have proved that the symmetric metric is
a critical point (in suitable sense) for the functional $\lambda_1$ on
the set of \emph{all} Riemannian metrics with fixed volume.

Curiously in the two exceptional examples (E-type) the best estimate
gotten by this method is strictly larger than 2, which is $\la_1$ of
the symmetric metric.
\begin{teo}\label{exceptional-intro}
If  $X= E_6/P(\alfa_1)$ 
(resp.   $X= E_7/P(\alfa_7)$) then
\begin{equation*}
  \la_1(X, g) \leq \frac{36}{17}
\qquad \text{ resp. } \quad 
  \la_1(X, g) \leq \frac{133}{53}
\end{equation*}
for any \Keler metric $g$ whose \Keler class $\om_g$ lies in
$2\pi\chern_1(X)$. 
\end{teo}
It would be interesting to understand if this is a deficiency of the
method, 
or if these symmetric spaces do in fact support metrics with $\la_1$
larger than 2.

\paragraph{Acknowledgements}
We wish to thank Prof. Gian Pietro Pirola and Prof. Andrea Loi for
various interesting discussions. The second author is also grateful to
Prof.  Peter Heinzner for inviting him to Ruhr--Universit\"at Bochum
during the preparation of this work. He also acknowledges partial
support by MIUR COFIN 2005 ``Spazi di moduli e teoria di Lie''.

\section{Stability of the Gieseker point}

Let $X$ be a compact complex manifold and let $E\ra X$ be a
holomorphic vector bundle of rank $r$.  Set
\begin{equation}
  V=H^0(X,E), \quad V'=H^0(X,\det E) \quad \w = \Hom(\Lambda^r V, V').
\end{equation}
The algebraic group $\Gl(V)$ acts linearly on $V$ hence on $\Lambda^r
V$.  It therefore also acts on $\w$.  For
$a\in \Gl(V)$ let $\pot[a] $ be the induced map on $\Lambda^r V$.  The
action of $\Gl(V)$ on $\w$ is given by
\begin{equation}
  a.T := T \comp (\pot[a])\meno.
\label{eq:GLact}
\end{equation}
Consider the above action restricted to the subgroup $\Sl(V) \subset
\Gl(V)$.  According to the terminology of geometric invariant theory,
a point $T\in \w $ is \enf{stable} (for this restricted action) if the
orbit of $\Sl(V)$ through $T$ is closed in $\w$ and the stabiliser of
$T$ inside $\Sl(V)$ is finite.  We denote by $\ts$ and $S$ the
stabilisers of $T_E$ in $\Gl(V)$ and $\Sl(V)$ respectively:
\begin{equation}
  \label{eq:def-S}
  \ts=\{a\in \Gl(V): a.T_E=T_E\} \qquad  S=\{a\in \Sl(V): a.T_E=T_E\}.
\end{equation}
For $x\in X$ put
\begin{gather*}
  V_x = \{ s\in V: s(x) = 0\} \qquad V'_x = \{ t\in V': t(x) = 0\} .
\end{gather*}
The following two simple lemmata will be used in the (algebraic) proof
of Thm. 1.
\begin{lemma} \label{genero} 
  Let $E$ be globally generated of rank $r$. Then a section $s\in V$
  belongs to $V_x$ if and only if for any choice of $r-1$ sections $
  s_2, \ldots{}, s_r$ of $E$ the section $ T_E(s, s_2, \ldots{} , s_r
  )$ of $\det E$ lies in $ V'_x$.
\end{lemma}
The proof is immediate.
Let $\Aut(E)$ be the group of holomorphic bundle automorphisms of $E$.
If $f\in \Aut(E)$ and $x\in X$, the map $f_x : E_x \ra E_x$ is a
linear isomorphism. The function $x\mapsto \det f_x$ is holomorphic
hence a (nonzero) constant and $f\mapsto \det f$ is a character of
$\Aut(E)$.  For $f\in \Aut(E)$ and $s \in H^0(X,E)$ put
\begin{equation}
  \eps(f) (s)  := f \comp s.
\end{equation}
This defines a representation $\eps : \Aut(E) \ra \Gl(V)$.
\begin{lemma} \label{lemma-aut} Let $E$ be globally generated.  Then $
  \eps( \{ f\in \Aut(E) : \det f =1\}) =\ts$.
\end{lemma}
\dimo For $f\in \Aut(E)$
\begin{equation*}
  \begin{gathered}
    \bigl ( \eps(f). T_E\bigr ) (s_1 , \ldots{}, s_r )(x) =
    \bigl( \eps(f)\meno s_1\bigr ) (x) \wedge \cdots{} \wedge \bigl(
    \eps(f)\meno s_r\bigr ) (x)  = \\
    = \bigl ( f_x\meno s_1(x) \bigr ) \wedge \cdots{} \wedge \bigl (
    f_x\meno s_r(x) \bigr )
    = \bigl ( \det f \bigr )\meno \cdot s_1(x) \wedge \ldots{} \wedge
    s_r(x) = \\
    =\bigl ( \det f \bigr )\meno \cdot T_E(s_1 , \ldots{}, s_r)(x).
  \end{gathered}
\end{equation*}
So
\begin{equation}
  \label{eq:computa-computa}
  \eps(f). T_E = \bigl ( \det f \bigr )\meno \cdot T_E.
\end{equation}
If $\det f=1$ then $\eps(f).T_E =T_E$. This proves that $ \eps( \{
f\in \Aut(E) : \det f =1\}) \subset \ts$.  Conversely, let $a\in \ts$.
We claim that $a ( V_x) = V_x$ for any $x\in X$.  Indeed let $s\in
V_x$. Then for any $s_2, \ldots{}, s_r \in V$
\begin{gather*}
  T_E(as, s_2, \ldots{}, s_r) = T_E(as, aa\meno s_2 , \ldots{},
  aa\meno s_r) =\\
  =(a\meno .\, T_E) (s, a\meno s_2, \ldots{}, a\meno s_r)=\\
  = T_E(s, a\meno s_2, \ldots{}.  a\meno s_r)(x) =0.
\end{gather*}
So $a s \in V_x$ by Lemma \ref{genero}, and indeed $a(V_x) =V_x$ as
claimed. We get therefore an induced isomorphism
\begin{gather*}
  f_x : E_x \cong V/ V_x \ra V/ V_x \cong E_x.
\end{gather*}
By construction $f_x ( s(x)) = (a\, s)(x)$. Since $E$ is globally
generated this ensures that $f$ is holomorphic so $f \in \Aut(E)$ and
$\eps(f)=a$. By \eqref{eq:computa-computa} $ \bigl ( \det f \bigr
)\meno \cdot T_E = \eps(f) .T_E = a.T_E =T_E$. Since $E$ is globally
generated, $T_E \neq 0$ and it follows that $\det f =1$. \fine
We recall two results that will be needed in the following.
\begin{teo}\label{luna}
  Let $H$ be a reductive group and $K \subset H$ a reductive subgroup.
  Let $X$ be an affine $H$-variety. If $x\in X$ is a fixed point of
  $K$ the orbit $Hx$ is closed if and only if the orbit $N_H(K)x$ is
  closed.
\end{teo}
This criterion is due to Luna \cite[Cor. 1]{luna-inv-75} and is based
on the Slice Theorem. For a complex analytic proof of the Slice
Theorem see \cite{heinzner-schwarz-Cartan}.

We recall that a \enf{rational homogeneous space} is a projective
variety $X$ of the form $G/P$ with $G$ a simply connected complex
semisimple Lie group and $P$ a parabolic subgroup without simple
factors. (See for example \cite {akhiezer-libro},
\cite{ottaviani-rat}, \cite{baston-eastwood}.)  Such spaces are also
called \enf {generalised flag manifolds}. A homogeneous vector bundle
$E$ over $X$ is of the form $E= G \times_P \V$ where $\V$ is a
representation of $P$.  If the representation is irreducible the
vector bundle itself is called \enf{irreducible}.

\begin{teo}[\protect {\cite[Thm.  1]{ramanan-homogeneous}}]
  \label{ramajan} Let $E\ra X$ be an irreducible homogeneous vector
  bundle over a rational homogeneous space $X$. Then $E$ is simple,
  i.e.  $\Aut(E) = \C^* \cdot I_E$.
\end{teo}
\noindent\textbf{First proof of Thm.   \ref{gieseker-intro}.}
Let $X=G/P$ be as above.  By Bott-Borel-Weil theorem (Thm.  \ref{bbw}
below) the hypothesis $H^0(E) \neq 0$ already ensures that $E$ is
globally generated, so both $V$ and $V'$ have positive dimension and
$T_E \not \equiv 0$.  By the same theorem $G$ acts irreducibly on both
$V$ and $V'$.  Denote by $ \rho : G \ra \Gl(V) $ and $ \sigma : G\ra
\Gl(V') $ these representations.  Since $G$ is semisimple
all the characters of $G$ are trivial.  In particular any
representation of $G$ on a vector space $U$ has
image contained in $\Sl(U)$.  So in fact $ \rho : G \ra \Sl(V) $ and $
\sigma : G\ra \Sl(V')$.  The Gieseker point is $G$-equivariant, that
is
\begin{equation} \label{equivario} T_E( \pot\rho(g) (u )) = \sigma(g)
 (  T(u) ) \qquad u \in \Lambda^r V.
\end{equation}
Set $H= \Sl(V) \times G$ and define a representation $ \varpi : H \ra
\Gl(\W) $ by
\begin{equation*}
  \varpi(a,g) T = \sigma(g)\comp T \comp \pot[a]\meno. 
\end{equation*}
Let $\tau : G \ra H$ be the morphism $\tau(g)=(\rho(g), g)$ and let
$\gh$ be the image of $\tau$.  $\gh\subset H$ is a closed reductive
subgroup and by \eqref{equivario} $T_E$ is a fixed point of $\gh$
acting via $\varpi$.  We claim that the normaliser $N_H(\gh)$ is a
finite extension of $\gh$.  In fact, denote by $\Ad$ the conjugation
on $H$.  Given $n\in N_H(\gh)$ put
\begin{equation}
  \phi(n) =\Ad(n)\restr{\gh} : \gh \ra \gh.
\end{equation}
Let $\Aut(\gh)$ denote the group of automorphisms of $\gh$ and
$\Inn(\gh)$ the subgroup of inner automorphisms.  Then $\phi: N_H(\gh)
\ra \Aut(\gh)$ is a morphism of groups.  Since $\gh$ is semisimple
$\Aut(\gh)$ is a finite extension of $\Inn(\gh)$ (see \cite [p.423]
{helgason} or
\cite[Thm. 1 p.203]{onishchik-vinberg-seminario}).  Put
$N'=\phi\meno(\Inn(\gh))$. Then $N' \lhd N_H(\gh)$ and
$$
N_H(\gh)/N' \hookrightarrow \Aut(\gh)/\Inn(\gh).
$$ 
Therefore $N_H(\gh)$ is a finite extension of $N'$ and it is enough to
prove that $N'$ is a finite extension of $\gh$.  Indeed if $n\in N'$
there is some $k\in K$ such that $nk'n\meno = kk'k\meno$ for any
$k'\in K$.  So $ k\meno n $ centralises $\gh$. If $k\meno n =(a,g) $
(with $a\in \Sl(V)$ and $g\in G$) this means that for any $g'\in G$ we
have
\begin{gather*}
  a\rho(g') = \rho(g') a \qquad
  \qquad gg'=g'g.
\end{gather*}
The second formula says that $g\in Z(G)$.  The first formula says that
$a: V \ra V$ commutes with the representation $\rho $.  Since this is
irreducible Schur lemma implies that $a=\eps I$ for some $\eps \in
\C^*$.  But $a\in \Sl(V)$, so $\eps^p=1$ where $p=\dim V$. Denote by
$U_p$ the group of $p$--roots of unity. Then $k\meno n= (\eps,g) \in
U_p \times Z(G)$. This proves that the composition
\begin{equation*}
  U_p\times Z(G) \rightarrow N' \ra N'/ \gh
\end{equation*}
is onto. Since $Z(G)$ is finite, it follows that $N'$ and $N_H(K)$ are
finite extensions of $K$.  Now $T_E\in \W$ is a fixed point of $\gh$
and $N_H(\gh)$ is a finite extension of $\gh$, so the orbit $N_H(\gh).
T_E$ is a finite set, hence it is closed.  Notice that both $H$ and
$K$ are reductive.  We can therefore apply Luna's criterion (Thm.
\ref{luna}) to the effect that the orbit $H.T_E$ is closed.  To finish
we claim that $H.T_E = \Sl(V).T_E$.  Since the action \eqref{eq:GLact}
of $\Sl(V)$ and the restriction of $\varpi$ to $\Sl(V)\times \{1\}
\subset H$ agree, the inclusion $H.T_E \supset \Sl(V).T_E$ is obvious.
For the other let $h=(a,g)\in H$. Then
\begin{gather*}
  \varpi(h) T_E = \sigma(g) \comp T_E \comp \pot[a]\meno = \sigma(g)
  \comp T_E \comp   \pot\bigl (a \rho(g)\meno \cdot \rho( g) \bigr)\meno =\\
  = \sigma(g) \comp T_E \comp \pot\rho(g)\meno \comp \pot ( \rho(g)
  a\meno) = \\
  =T_E \comp \pot ( \rho(g) a\meno) = \varpi(a\rho(g)\meno, 1)T_E
\end{gather*}
and $a\rho(g\meno) \in \Sl(V)$. Therefore $H.T_E \subset \Sl(V).T_E$
so the two orbits coincide.  This shows that the orbit of $T_E$ is
closed. Let $S$ and $\ts$ be the stabilisers defined as in
\eqref{eq:def-S}.  By Thm. \ref{ramajan}, $\Aut(E) = \C^* \cdot I_E$,
therefore $\{f\in \Aut(E): \det f=1\}$ is finite, which implies, by
Lemma \ref{lemma-aut}, that $\ts$ and a fortiori $S$ are finite.
\fine We remark that this proof works over any algebraically closed
field of characteristic zero.

We come now to the second proof of this result.  Recall that if $E$ is
a globally generated bundle on $X$ and $\boldsymbol{s}=\{s_1,
\ldots{}, s_N\}$ is a basis of $H^0(X,E)$ there is an induced map $
\phi_{\boldsymbol{s}} : X \ra G(r,N)$.  Consider on $G(r,N)$ the
standard symmetric K\"ahler structure which coincides with the
pull-back of the Fubini-Study metric via the Pl\"ucker embedding.
Denote by $\mu: G(r,N) \ra \su(N)$ the moment map for the standard
action of $\SU(N) $ on $G(r,N)$.
\begin{teo} [{\cite[Thm.  3.1]{wang-xiaowei-balance}}]
  \label{wang}
  Let $(X^m,\om)$ be a compact \Keler manifold and let $E$ be a
  globally generated bundle on $X$.  Then $T_E$ is stable if and only
  if there is a basis $ {\boldsymbol{s}} $ of $H^0(X,E)$ such that
  \begin{equation}
\label{omega-bil}
    \int_X \mu \bigl ( \phi_  {\boldsymbol{s}} (x) \bigr ) \om^m (x) =0.    
  \end{equation}
\end{teo}
For the reader's convenience we briefly sketch the proof. \\
\noindent
{\bf  Proof.\ }
Fix an arbitrary Hermitian metric $h$ on $E$ and
consider on $V$ the corresponding 
 $L^2$--scalar product. 
Let
$\basis{s}$ be an orthonormal basis with respect to this product.  On
the line bundle $\det E$ consider the metric $\phi_\basis{s}^* h_G$
where $h_G$ is the metric on $\mathcal{O}_{G(r,N)}(1)$.  
Consider on $V'$ the corresponding  
$L^2$--scalar product. 
Finally denote by  $\scalar_{\W}$
the Hermitian inner product on $\W$,
 $|| \cdot ||_\w$ being the corresponding norm.  
 Since we have fixed a basis we may identify $\Sl(V)$ with $\Sl(N,
 \C)$.  For $g\in \Sl(N, \C)$ set $\nu(g) = \log || g\meno . T_E
 ||_\W$. We consider $\nu$ as a function on $\Sl(N, \C)/ \SU(N)$.  On
 this space Wang introduces another functional
\begin{equation*}
  L(g):= \int_M \biggl ( \sum_I ||  (g^{-1}T_E) (s_{I})(x)
  ||^2_{\phi_\basis{s}^*h_G} \bigr )
  \frac{\om^n}{n!}(x),
\end{equation*}
which is strictly convex on $\Sl(N,\C)/SU(N)$ \cite[Lemma
3.5]{wang-xiaowei-balance}.  (Here $s_I=s_{i_1} \wedge \cdots \wedge s_{i_r}\in \Lambda^r V$.)
Critical points of $L$ correspond to $g\in \Sl(N)$ such that the basis
$\{gs_1, \ldots{}, gs_N\} $ satisfies \eqref{omega-bil}.  For some
constants $C_2, C_4 \in \mathbb{R}$ and $C_1, C_3 >0$ the inequalities
\begin{equation}
L
 \geq C_1 \nu + C_2  \geq C_3 L + C_4
\label{wang-puppa}
\end{equation}
hold on $ \Sl(N,\C)/SU(N)$.  The first is proved by Wang \cite[
p.406]{wang-xiaowei-balance}.  The second is simply an application of
Jensen inequality to the convex function $-\log$.  If $T_E$ is stable,
then $\nu$ is proper by the Kempf-Ness theorem \cite {kempf-ness}.
Hence $L$ is proper too, so admits a minimum and there is a basis
$\basis{s'}$ such that \eqref{omega-bil} is satisfied. On the other
hand if there is such a basis, $L$ has a minimum and being strictly
convex this means it is proper. By \eqref{wang-puppa}, $\nu$ is proper
as well and, again by Kempf-Ness theorem, this implies that $T_E$ is
stable.  It should be noted that the identification of the moment map
for a projective action with the differential of a convex functional
is standard in analytic Geometric Invariant Theory
\cite[Ch. 8]{mumford-GIT}, \cite[\S 6.5]{donaldson-kronheimer},
\cite{heinzner-huckleberry-MSRI}.  \fine

\noindent\textbf{Second proof of Theorem \ref{gieseker-intro}.}
Let 
$K$ be a compact form of $G=\Aut(X)$.  By averaging on $K$ we can find
$K$-invariant metrics $\om$ and $h$ on $X$ and $E$ respectively. Let
$\boldsymbol{s}$ be a basis of $H^0(X,E)$ that is orthonormal with
respect to the $L^2$-scalar product obtained using $h$ and $\om$. By
Bott-Borel-Weil theorem (Thm. \ref{bbw} below) $G$ and hence $K$ act
irreducibly on $H^0(X,E) \cong \C^N$. Denote by $\sigma : K \ra
\SU(N)$ this representation (recall that $K$ is semisimple).  Then
$\mu \comp \phi_{\boldsymbol{s}} $ is $K$-equivariant and
\begin{equation*}
  B=\int_X \mu \bigl ( \phi_  {\boldsymbol{s}} (x) \bigr ) \om^m (x) 
\end{equation*}
is a fixed point of $\ad(\sigma(K)) \subset \Gl( \su(N))$, that is
$\sigma(k) B = B \sigma(k)$ for $k\in K$. By Schur lemma this implies
that $B=\la I$, so $B=0$ since $B\in \su(N)$. By Thm. \ref{wang} the
Gieseker point $T_E$ is stable. \fine


In order to clarify the meaning of the above result it might be good
to notice that together with the numerical criteria of
\cite{gieseker-vector-surfaces} it allows an easy proof of the
Gieseker stability (see e.g.  \cite[p.189]{kobayashi-vector}) of
irreducible homogeneous bundles. We sketch this argument, although a
stronger result (Mumford-Takemoto stability) is well-known (see
\cite[p.65]{ottaviani-rat} and references therein).

.
\begin{prop}
  [{\cite[Prop. 2.3]{gieseker-vector-surfaces}}] \label{gieseker123}
  Let $T\in \w$ be a stable point.  Let $V''\subset V$ be a subspace
  and let $d$ a number $1 \leq d < r $.  Assume that for any $d+1$
  vectors $v_1, \ldots{}, v_{d+1} \in V''$, $T(v_1, \ldots{}, v_d,
  v_{d+1} , \cdots) \equiv 0$. Then
 $   \dim V'' < (d/r) \, \cdot \,  \dim V.$
  If $T$ is only semistable, then equality can hold.
\end{prop}
In \cite{gieseker-vector-surfaces} there is a proof in the semistable
case, which works as well in the stable case.
\begin{cor}
  Let $E\ra X$ be an irreducible homogeneous vector bundle of rank $r$
  over a rational homogeneous space $X=G/P$.  If $H^0(E) \neq 0$, and
  $F \subset E$ is a subsheaf of rank $d$, then
  $
    h^0(F) < (d/r) \cdot h^0(E).
  $
\end{cor}
Fix now an irreducible homogeneous bundle $E$ of rank $r$ and let
$F\subset E$ be a subsheaf of rank $d$, with $0<d<r$. Let $H$ be any
polarisation on $X$. Since any line bundle is homogeneous, $E(k)=E
\otimes H^{\otimes k}$ is homogeneous.  By Serre Theorem there is a
$k_0$ such that for $k\geq k_0$
$$
H^i(X, F(k) ) = H^i(X, E(k))=\{0\} \qquad i > 0
$$
and both $E(k)$ and $F(k)$ are globally generated. By
Thm. \ref{gieseker-intro}, $T_{E(k)}$ is stable, so by the above
corollary $ \chi(X, F(k)) = h^0(X,F(k)) < (d/r) \cdot h^0(X,E(k)) =
(d/r) \cdot \chi(X,E(k)).  $ This proves that any irreducible
homogeneous bundle is Gieseker stable with respect to any
polarisation.
%
%
%

\section{The first eigenvalue of Hermitian symmetric \\
  spaces}

Here we want to apply the previous stability result to a problem in
spectral geometry. Let $X$ be a projective manifold and $L$ an ample
line bundle on $X$.  Let $\m$ be the set of \Keler metrics $g$ with
\Keler form $\om_g$ lying in the class $2\pi\chern_1(L)$.  For $g$ in
$\m$ let $\Delta_g $ be the Laplacian on functions,
\begin{equation*}
  \Delta_g f = -d^* d f = 
  2 \  g^{i\bar{j}}\frac{\partial^2 f}{\partial z^i \partial \bar{z}^j}.
\end{equation*}
It is well-known that $\Delta_g$ is a negative definite elliptic
operator and has therefore discrete spectrum: denote its eigenvalues
by $0 > -\la_1(g) > -\la_2(g) > \cdots$.  The following result of
Lichnerowicz relates $\la_1$ to \KE geometry.
\begin{teo}[{\cite[Thm.  2.4.3, p.41]{futaki-libro}}] \label{futaki}
  If $X$ is a Fano manifold and $g_{KE}$ is a \KE metric, i.e.
  $Ric(g_{KE})=g_{KE}$, then $\la_1(g_{KE}) =2$ if $\Aut(X)$ has
  positive dimension and $\la_1(g_{KE}) > 2$ otherwise.
\end{teo}
We are interested in upper estimates for $\la_1(g)$ of general metrics
in the class $\m$.  Bourguignon, Li and Yau \cite{bourguignon-li-yau}
first studied this problem and showed that the supremum
\begin{equation}\label{eq:def-sup}
  I(L)=  \sup_{\m} \la_1(g) 
\end{equation}
is always finite. (This heavily depends on the restriction to \Keler
metrics, see \cite{colbois-dodziuk}.)  They gave an explicit upper
bound for $I(L)$ in terms of numerical invariants of a globally
generated line bundle $E$. For $(X,L)=(\mathbb{P}^m,
\mathcal{O}_{\mathbb{P}^m} (1))$ they were able to show that $I(L)=2$.
The following criterion, due to Arezzo, Loi and the second author, is
an extension of Bourguignon, Li and Yau's theorem. It allows to attack
this problem using holomorphic vector bundles instead of just line
bundles.
\begin{teo}[{\cite[Thm. 1.1]{arezzo-ghigi-loi}}]\label{igi}
  Let $(X,L)$ be a polarised manifold and $E $ a holomorphic vector
  bundle of rank $r$ over $X$. 
  Assume that $E$ is globally generated and nontrivial and put
  \begin{equation}
    \label{eq:def-J}
    J(E,L):=    \frac{2\, \dim_\C  X \cdot h^0(E)\,  
      \langle c_1(E) \cup c_1(L)^{m-1} ,[X]\rangle  }
    {r\, (h^0(E)-r)\, \langle c_1(L)^m, [X]\rangle}.
  \end{equation}  
  If the Gieseker point $T_E$ is stable, then
  \begin{equation}
    \label{mainest1}
    I(L)\leq J(E,L).
  \end{equation}
\end{teo}
The result of \cite{arezzo-ghigi-loi} is slightly more general since
there is no projectivity assumption on $X$.

We want to apply this result to the case where $X=G/P$ is a rational
homogeneous space and $E$ is homogeneous. In this case
$J(E,L)$ 
can be computed, at least in principle, in terms of Lie algebra data.
To proceed we fix the following (standard) notation. (See e.g.
\cite[Ch.  1]{fels-huckleberry-wolf},\cite{ottaviani-rat},
\cite{baston-eastwood}.)  $G$ is a simply connected complex semisimple
Lie group, $\lieg = {Lie}\, G$, $\lieh \subset \lieg$ is a Cartan
subalgebra, $l=\dim \lieh$ is the rank of $G$, $B$ is the Killing form
of $\lieg$, $\Delta$ is the root system of $(\lieg, \lieh)$,
$\Delta_+$ is a system of positive roots, $\Delta_-=-\Delta_+$,
$\Pi=\{\alfa_1, \ldots{}, \alfa_l\}$ is the set of simple roots,
$\fund_1, \ldots, \fund_l$ denote the fundamental weights.  $\Lambda=
\Zeta \fund_1 \oplus \cdots \oplus \Zeta \fund_l \subset \lieh^*$ is
the weight lattice of $\lieg$ relative to the Cartan subalgebra
$\lieh$.  For $\alfa \in \Delta$ let $H_\alfa \in \lieh$ be such that
$\alfa(X) = B(X,H_\alfa)$.  $\lieb$ is the standard \enf{negative}
Borel subalgebra:
\begin{equation}
  \label{eq:borel}
  \lieb=\lieh \oplus \bigoplus _{\alfa \in \Delta_-} \lieg_\alfa
\end{equation}
Parabolic subalgebras containing $\lieb$ are of the form
\begin{equation}
  \liep(\Sigma) = \lieb \oplus \bigoplus _{\alfa \in \spam (\Pi - \Sigma) \cap
    \Delta_+} \lieg_\alfa
\end{equation}
where $\Sigma $ is some subset of $ \Pi$. For example $\Sigma = \Pi$
corresponds to $\lieb$, $\Sigma = \varnothing$ to $\lieg$ and
maximally parabolic subalgebras are of the form $ \liep(\alfa_k)$.
The algebra $\liep(\Sigma)$ admits a Levi decomposition $
\liep(\Sigma) = \liel(\Sigma) \oplus \mathfrak{u}(\Sigma) $, where
$\mathfrak{u}(\Sigma)$ is the nilpotent radical and
$$
\liel(\Sigma) = \lieh \oplus \bigoplus_{\alfa \in \spam(\Pi-\Sigma)
  \cap \Delta} \lieg_\alfa
$$
is the reductive part. This latter admits a further decomposition $
\liel(\Sigma ) = \liez(\Sigma) \oplus \lies(\Sigma)$, $\liez(\Sigma) $
being the center and $\lies(\Sigma)$ being semisimple. Moreover
\begin{equation*}
  \liez(\Sigma) = \bigcap_{\alfa \in\Pi-\Sigma} \ker \alfa \subset \lieh
\end{equation*}
and
\begin{equation*}
  \lies(\Sigma) = \lieh' (\Sigma) \oplus \bigoplus_{\alfa \in \spam(\Pi-\Sigma) \cap
    \Delta} \lieg_\alfa
\end{equation*}
where $\lieh'(\Sigma)= \spam \{ H_\alfa : \alfa \in \Pi- \Sigma \}
\subset \lieh $ is a Cartan subalgebra for $\lies(\Sigma)$ and
$\lieh=\liez(\Sigma) \oplus \lieh'(\Sigma)$.  We denote by $B,
P(\Sigma), L(\Sigma), U(\Sigma), Z(\sigma), S(\Sigma)$ the
corresponding closed subgroups of $G$. Note that $S(\Sigma)$ is simply
connected.  One can describe $\liep(\Sigma), P(\Sigma)$ and the
homogeneous space $G/P(\Sigma)$ by the Dynkin diagram of $G$ with the
nodes corresponding to roots in $\Sigma$ crossed.

A weight $\la=\sum_i m_i \fund_i \in \Lam $ is \enf{dominant for $G$}
or simply \enf{dominant} if $m_i \geq 0$ for any $i$. It is said to be
\enf{dominant with respect to $\liep(\Sigma)$} if $m_i\geq 0$ for any
index $i$ such that $\alfa_i\not \in \Sigma$. By highest weight
theory, the irreducible representations of $G$ are parametrised by
dominant weights, while irreducible representations of a parabolic
subgroup $P(\Sigma)$ are parametrised by weights that are dominant
with respect to $\liep(\Sigma)$.  If $\la $ is dominant we let $W_\la$
denote the irreducible representation of $G$ with highest weight
$\la$.  If $\la $ is dominant for $\liep(\Sigma)$ we let $V_\la$
denote the irreducible representation of $P(\Sigma)$ with highest
weight $\la$. We let moreover $E_\la$ denote the homogeneous vector
bundle on $X=G/P(\Sigma)$ defined by the representation $V_\la$, that
is $ E_\la = G \times_{P(\Sigma)} V_\la.  $
\begin{teo}
  [Bott-Borel-Weil] \label{bbw} If $\la \in \Lam$ is dominant for $G$,
  then
$$
 H^0(X, E_\la) = W_\la.
$$
Otherwise $H^0(X,E_\la)=\{0\}$.
\end{teo}
Bott's version of the theorem is much more general, but this partial
statement is enough for what follows.  We also remark that if one
chooses $\lieb$ to be the Borel subalgebra with \enf{positive} instead
of negative roots, which is customary for example in the usual picture
of $\mathbb{P}^n$ as the set of lines in $\C^{n+1}$, then one has to
consider \enf{lowest weights} instead of highest ones. This amounts to
dualize both representations.  With this choice the statement of the
theorem becomes $ H^0(X, E^*_\la) = \bigl ( W_\la\bigr )^*.$ (The book
\cite{baston-eastwood} follows this convention.)

Recall that the set of simple roots $\Pi=\{\alfa_1, \ldots, \alfa_l\}$
is a basis of $\Lam\otimes \mathbb{Q}$. For a weight $\la\in \Lam$,
let $\la=\sum_i \xi_i (\la)\alfa_i$ be its expression in this basis.
We say that the (rational) number $\xi_i(\la)$ is the \enf{coefficient
  of $\alfa_i$ in $\la$}.
We denote by $\la_\ad$ the highest weight of the adjoint
representation of $G$ (that is the largest root).
\begin{lemma} \label{c1} Let $X=G/P(\alfa_k)$. The bundle
  $E_{\fund_k}$ associated to the fundamental weight $\fund_k$ is a
  very ample line bundle over $X$.  Moreover
  $\operatorname{Pic}(X)\cong H^2(X,\Zeta)=\Zeta
  \chern_1(E_{\fund_k})$. For any weight $\la \in \Lam$ that is
  dominant for $P(\alfa_k)$
  \begin{equation}
    \chern_1(E_\la) = \dim V_\la \, \frac{\xi_k(\la)}{\xi_k(\fund_k)} \, \chern_1(E_{\fund_k}).
  \end{equation}
\end{lemma}
(For the proof see e.g. \cite[\S 5.2]{ramanan-homogeneous}, \cite[p.56]{ottaviani-rat}.)  

In the following statement we summarise what we need of the structure
theory of Hermitian symmetric spaces.
\begin{teo}\label{classifico}
  An irreducible Hermitian (globally) symmetric space of the compact
  type is a rational homogeneous space. Moreover a rational
  homogeneous space $X=G/P$ is symmetric if and only if the
  representation of $\liep$ on $\lieg / \liep$ induced from the
  adjoint representation of $\lieg$ is irreducible. The actual
  possibilities are explicitly listed in Table 1.
\end{teo}
The characterisation in terms of irreducibility of $\lieg/\liep$ is
due to Kobayashi and Nagano \cite[Thm.
A]{kobayashi-nagano-filtered-II} (see also
\cite[p.26]{baston-eastwood}).

\begin{table}[h]
  \centering
  \label{tab1}
  \begin{tabular}[c]{|p{0.2\linewidth}|c|c|c|}
    \hline{}   
    & Klein form & & Type \\
    \hline{}
    Grassmannian & $G_{k,n}= \Sl(n)/P(\alfa_k)$ & $n\geq 2$ & AIII \\
    \hline{} 
    Odd quadrics & $Q_{2n-1} = \Spin(2n+1) /P(\alfa_1)$ & $n\geq 2 $ & BI \\
    \hline{} 
    Even quadrics & $Q_{2n-2} = \Spin(2n) /P(\alfa_1)$ & $n\geq 3 $ & DI \\
    \hline{} 
    Spinor variety & $X=\Spin(2n)/P(\alfa_n)$ & $n\geq 4$ & DIII\\
    \hline{}
    Lagrangian
    Grassmannian     
    &
    $X= \Sp(n,\C) / P(\alfa_n)$ 
    & $n \geq 2$ & CI \\
    \hline
    & $X= E_6 / P(\alfa_1)$ & & EIII \\
    \hline
    & $X= E_7 / P(\alfa_7) $ & & EVII \\
    \hline
  \end{tabular}
  \caption{Irreducible Hermitian symmetric spaces of the compact type.}
\end{table}

\begin{prop} \label{prop-J} Let $X=G/P(\alfa_k)$ be a compact
  irreducible Hermitian symmetric space and let $\la\in \Lam$ be a
  nontrivial dominant weight.  Then
  \begin{equation}
    \label{eq:J-hom}
    J(E_\la, -K_X) = \frac{ 2\, \dim W_\la } { \dim W_\la - \dim V_\la } \cdot
    \frac{\xi_k(\la)}{\xi_k(\la_\ad)}.
  \end{equation}
\end{prop}
\dimo The tangent bundle to $X=G/P$ is the homogeneous bundle obtained
from the representation of $P$ on $\lieg/\liep$.  For symmetric $X$
this is irreducible by Theorem \ref{classifico}, so Bott-Borel-Weil
theorem and Lemma \ref{c1} apply. Since $H^0(X,TX)= \lieg =
W_{\la_\ad} $ (see \cite[p.75, p.131]{akhiezer-libro}),
$\lieg/\liep = V_{{\la_\ad}}$ and $TX=E_{\la_\ad}$.  Set $m=\dim X =
\dim V_{\la_\ad} $.  By Lemma \ref{c1}
\begin{gather*}
  \chern_1(-K_X) =\chern_1(TX)= m \cdot
  \frac{\xi_k(\la_{\ad})}{\xi_k(\fund_k)} \, \chern_1(E_{\fund_k}) , \\
  \chern_1(E_\la) = \dim V_\la \cdot \frac{\xi_k(\la)}{\xi_k(\fund_k)}
  \, \chern_1(E_{\fund_k})  , \\
  \frac{ \langle c_1(E_\la) \cup c_1(-K_X)^{m-1} ,[X]\rangle } { \langle
    c_1(-K_X)^m, [X]\rangle} = \frac{\dim V_\la}{m} \cdot
  \frac{\xi_k(\la)}{\xi_k(\la_\ad)}.
\end{gather*}
The rank of $E_\la$ is $\dim V_\la$, while $h^0(X,E_\la)=\dim W_\la$
by Bott-Borel-Weil theorem. Therefore
\begin{gather*}
  J(E_\la,-K_X)= \frac{2\, m\, h^0(E_\la)} {r\, (h^0(E_\la)-r)} \cdot
  \frac{ \langle c_1(E_\la) \cup c_1(-K_X)^{m-1} ,[X]\rangle }
  { \langle c_1(-K_X)^m, [X]\rangle} = \\
  = \frac{ 2\, m\, \dim W_\la} { \dim V_\la ( \dim W_\la - \dim V_\la)
  } \cdot \frac{\dim V_\la}{m} \cdot
  \frac{\xi_k(\la)}{\xi_k(\la_\ad)} = \\
  = \frac{ 2\, \dim W_\la } { \dim W_\la - \dim V_\la } \cdot
  \frac{\xi_k(\la)}{\xi_k(\la_\ad)}.
\end{gather*}
\fine
We are now ready for the proof of theorems 2 and 3. \\
\noindent\textbf{Proof of Theorem \ref{classical-intro}.}
Let $X$ be a compact irreducible Hermitian symmetric space. Denote by
$g_{KE}$ the symmetric (K\"{a}hler-Einstein) metric with \Keler form
in $2\pi c_1(X)$.  We need to show that
\begin{equation}
  I(-K_X) = 2 = \la_1(g_{KE}).
\end{equation}
The second equality follows from Thm. \ref{futaki}. So $ I(-K_X) \geq
2$ by definition \eqref{eq:def-sup}. It is enough to prove that
$I(-K_X) \leq 2$.  For each space in the first five families in Table
\ref{tab1} we find a homogeneous bundle $E_\la\ra X$ such that
$J(E_\la, -K_X)=2$.  The result then follows applying Thm.  \ref{igi}.
The relevant information regarding weights and degrees can be found
for example in \cite[p.66, p.69]{humphreys-algebras}.
\\
\noindent{}
1. The case of the Grassmannians (type $AIII$) is settled by hand in
\cite[Thm. 1.3]{arezzo-ghigi-loi}. The vector bundle $E$ is the dual
of the universal subbundle. If we choose the Borel group as in
\eqref{eq:borel} then $ E=E_{\fund_1}$.
\\
\noindent{}
2. For odd quadrics the Dynkin diagram is:
\begin{center}
  \parbox[b]{0.45\textwidth}{
    $Q_{2n-1} = \Spin(2n+1) /P(\alfa_1)$ \\
    Type BI}
  \makebox[0.45\textwidth][b]{ \includegraphics{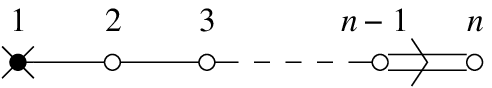}}
\end{center}
The largest root is $\la_\ad = \fund_2$.  Put $\la=\fund_n$.  Then
$W_\la$ is the spin representation, while $V_\la$ corresponds to the
spin representation of the semisimple part $S(\alfa_1) \cong
\Spin(2n-1)$ of $P(\alfa_1)$ . The bundle $E_\la$ is the \enf{spinor
  bundle} studied e.g. by Ottaviani \cite{ottaviani-spinor}.  Of
course $\dim W_\la = 2^{n}$, $\dim V_\la = 2^{n-1}$. Finally
$\xi_1(\la ) =\xi_1(\fund_n) = 1/2$, $\xi_1(\la_\ad ) =\xi_1(\fund_2)
= 1$, so
\begin{equation*}
  J(E_\la, -K_X) = \frac{ 2\, \dim W_\la } { \dim W_\la - \dim V_\la } \cdot
  \frac{\xi_1(\la)}{\xi_1(\la_\ad)} = 
  2 \cdot \frac{2^{n}} {2^{n} -2^{n-1}} \cdot \frac{1/2}{1} = 2.
\end{equation*}
\\
\noindent{}
3. The situation is very similar for even quadrics. The Dynkin diagram
is:
\begin{center}
  \parbox[b]{0.45\textwidth}{
    $Q_{2n-2} = \Spin(2n) /P(\alfa_1)$ \\
    Type DI\\
    \, \\
    \, }
  \makebox[0.45\textwidth][b]{ \includegraphics{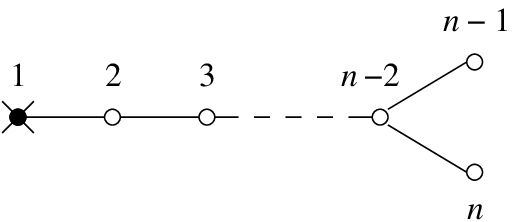}}
\end{center}
The largest root is again $\la_\ad = \fund_2$.  We take $W_\la$ to be
either one of the half-spin representation. ($E_\la$ is one of the two
spinor bundles on $Q_{2n-2}$, \cite{ottaviani-spinor}).  Say
$W_\la=\mathcal{S}_+$.  Then $\la = \fund_n$ and $V_\la$ is the
half-spin representation $\mathcal{S}_+$ of $S(\alfa_1) \cong
\Spin(2n-2)$.  Now $ \dim W_\la = 2^{n-1}$, $V_\la = 2^{n-2}$,
$\xi_1(\fund_n) = 1/2$, $\xi_1(\fund_2) = 1$, so again $J(E_\la, -K_X)
=2$.
\\
\noindent{}
4.  For the Lagrangian Grassmannian the Dynkin diagram is:
\begin{center}
  \parbox[b]{0.45\textwidth}{
    $X= \Sp(n,\C) / P(\alfa_n)$  \\
    Type CI}
  \makebox[0.45\textwidth][b]{ \includegraphics{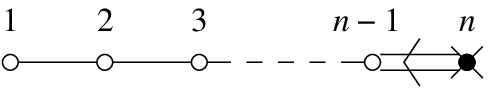}}
\end{center}
The highest weight of the adjoint representation is $\la_\ad = 2
\fund_1$. $W_{\fund_1}$ is the standard representation of $\Sp(n,\C)$
on $\C^{2n}$. The semisimple part of $P(\alfa_n)$ is
$S(\alfa_n)=\Sl(n)$, so $V_{\fund_1}$ is the standard representation
of $\Sl(n)$ on $\C^n$. So choosing $E=E_{\fund_1}$ we get
\begin{equation*}
  J(E, -K_X) = 2 \cdot \frac {2n}{2n -n} \cdot \frac{\xi_n(\fund_1) }
  {\xi_n(2\fund_1)} = 2.
\end{equation*}
\\
\noindent{}
5.  For the Spinor varieties the Dynkin diagram is:
\begin{center}
  \parbox[b]{0.45\textwidth}{
    $\Spin(2n)/P(\alfa_n)$ \\
    Type DIII\\
    \, \\ \,}
  \makebox[0.45\textwidth][b]{ \includegraphics{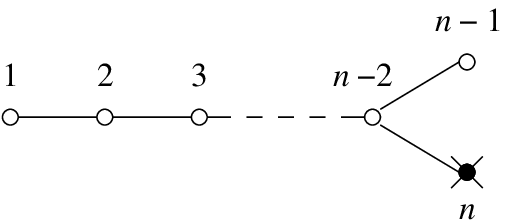}}
\end{center}
Take $E=E_{\fund_1}$. $W_{\fund_1}$ is the standard representation of
$\Spin(2n)$. The semisimple part of $P(\alfa_n)$ is
$S(\alfa_n)=\Sl(n)$, so $V_{\fund_1}$ is the standard representation
of $\Sl(n)$ on $\C^n$.  The largest root is $\la_\ad = \fund_2$,
$\xi_n(\fund_1)=1/2$, $\xi_n(\fund_2) = 1$. So
\begin{equation*}
  J(E, -K_X) = 
  2  \cdot \frac {2n}{2n -n} \cdot \frac{1/2 }{1} = 2.
\end{equation*}
\fine
\noindent\textbf{Proof of Theorem \ref{exceptional-intro}.}
1. For $X=E_6/P(\alfa_1)$ the Dynkin diagram (with Bourbaki numbering)
is:
\begin{center}
  \parbox[b]{0.45\textwidth}{ $ X=E_6/P(\alfa_1)$
    \\
    Type EIII\\
    \, }
  \makebox[0.45\textwidth][b]{ \includegraphics{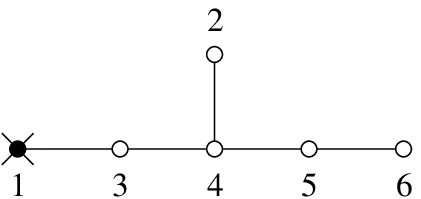}}
\end{center}
The largest root is $\la_\ad = \fund_2$.  An easy computation gives
$J(E_{\fund_6} , -K_X) $ $= 36/17$ and $J(E_{\fund_2} , -K_X) =
78/31$. If $\la = \sum_i a_i \fund_i$, then
\begin{gather*}
  J(E_\la, -K_X) \geq 2 \frac{\xi_1(\la)} {\xi_1(\fund_2)} =
  \frac{8}{3} a_1 + 2 a_2 + \frac{10}{3} a_3 + 4 a_4 + \frac{8}{3} a_5
  + \frac{4}{3} a_6 .
\end{gather*}
The right hand side is $< 36/17$ if and only if $\la=\fund_2$ or
$\la=\fund_6$. Therefore the best estimate is gotten with
$\la=\fund_6$.\\
\noindent{}
2. For $X=E_7/P(\alfa_7)$ the Dynkin diagram (with Bourbaki numbering)
is:
\begin{center}
  \parbox[b]{0.45\textwidth}{ $ X=E_7/P(\alfa_7)$
    \\
    Type EVII\\
    \, }
  \makebox[0.45\textwidth][b]{ \includegraphics{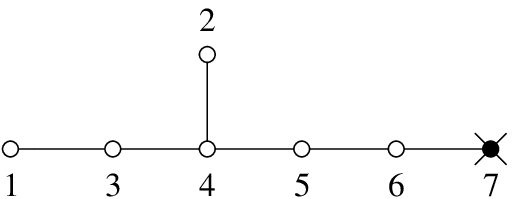}}
\end{center}
The largest root is $\la_\ad = \fund_1$.  We have $J(E_{\fund_1} ,
-K_X) = 133/53$. If $\la = \sum_i a_i \fund_i$, then
\begin{gather*}
  J(E_\la, -K_X) \geq 2 \frac{\xi_7(\la)} {\xi_7(\fund_1)} = 2
  \xi_7(\la)
  = \\
  = 2 a_1 + 3 a_2 +4 a_3 + 6 a_4 + 5 a_5 + 4 a_6 + 3 a_7.
\end{gather*}
The right hand side is $< 133/53$ if and only if $\la=\fund_1$.
Therefore the best estimate is gotten with $\la=\fund_1$.  \fine

\def\cprime{$'$}

\vskip.3cm

\noindent Universit\`a degli Studi di Parma,\\
\textit{E-mail:} \texttt{leonardo.biliotti@unipr.it}

\vskip.3cm

\noindent Universit\`a di Milano Bicocca,\\
\textit{E-mail:} \texttt{alessandro.ghigi@unimib.it}


\begin{thebibliography}{10}

\bibitem{akhiezer-libro}
D.~N. Akhiezer.
\newblock {\em Lie group actions in complex analysis}.
\newblock Aspects of Mathematics, E27. Friedr. Vieweg \& Sohn, Braunschweig,
  1995.

\bibitem{arezzo-ghigi-loi}
C.~Arezzo, A.~Ghigi, and A.~Loi.
\newblock Stable bundles and the first eigenvalue of the {L}aplacian.
\newblock {\em J. Geom. Anal.}, 17(3):375--386, 2007.

\bibitem{baston-eastwood}
R.~J. Baston and M.~G. Eastwood.
\newblock {\em The {P}enrose transform}.
\newblock Oxford Mathematical Monographs. The Clarendon Press Oxford University
  Press, New York, 1989.

\bibitem{bourguignon-li-yau}
J.-P. Bourguignon, P.~Li, and S.-T. Yau.
\newblock Upper bound for the first eigenvalue of algebraic submanifolds.
\newblock {\em Comment. Math. Helv.}, 69(2):199--207, 1994.

\bibitem{colbois-dodziuk}
B.~Colbois and J.~Dodziuk.
\newblock Riemannian metrics with large {$\lambda\sb 1$}.
\newblock {\em Proc. Amer. Math. Soc.}, 122(3):905--906, 1994.

\bibitem{donaldson-kronheimer}
S.~K. Donaldson and P.~B. Kronheimer.
\newblock {\em The geometry of four-manifolds}.
\newblock Oxford Mathematical Monographs. The Clarendon Press Oxford University
  Press, New York, 1990.
\newblock Oxford Science Publications.

\bibitem{el-soufi-ilias-Pacific}
A.~El~Soufi and S.~Ilias.
\newblock Riemannian manifolds admitting isometric immersions by their first
  eigenfunctions.
\newblock {\em Pacific J. Math.}, 195(1):91--99, 2000.

\bibitem{fels-huckleberry-wolf}
G.~Fels, A.~Huckleberry, and J.~A. Wolf.
\newblock {\em Cycle spaces of flag domains}, volume 245 of {\em Progress in
  Mathematics}.
\newblock Birkh\"auser Boston Inc., Boston, MA, 2006.
\newblock A complex geometric viewpoint.

\bibitem{futaki-libro}
A.~Futaki.
\newblock {\em K\"ahler-{E}instein metrics and integral invariants}.
\newblock Springer-Verlag, Berlin, 1988.

\bibitem{gieseker-vector-surfaces}
D.~Gieseker.
\newblock On the moduli of vector bundles on an algebraic surface.
\newblock {\em Ann. of Math. (2)}, 106(1):45--60, 1977.

\bibitem{heinzner-huckleberry-MSRI}
P.~Heinzner and A.~Huckleberry.
\newblock Analytic {H}ilbert quotients.
\newblock In {\em Several complex variables (Berkeley, CA, 1995--1996)},
  volume~37 of {\em Math. Sci. Res. Inst. Publ.}, pages 309--349. Cambridge
  Univ. Press, Cambridge, 1999.

\bibitem{heinzner-schwarz-Cartan}
P.~Heinzner and G.~W. Schwarz.
\newblock Cartan decomposition of the moment map.
\newblock {\em Math. Ann.}, 337(1):197--232, 2007.

\bibitem{helgason}
S.~Helgason.
\newblock {\em Differential geometry, {L}ie groups, and symmetric spaces},
  volume~80 of {\em Pure and Applied Mathematics}.
\newblock Academic Press Inc., New York, 1978.

\bibitem{humphreys-algebras}
J.~E. Humphreys.
\newblock {\em Introduction to {L}ie algebras and representation theory},
  volume~9 of {\em Graduate Texts in Mathematics}.
\newblock Springer-Verlag, New York, 1978.
\newblock Second printing, revised.

\bibitem{kempf-ness}
G.~Kempf and L.~Ness.
\newblock The length of vectors in representation spaces.
\newblock In {\em Algebraic geometry (Proc. Summer Meeting, Univ. Copenhagen,
  Copenhagen, 1978)}, volume 732 of {\em Lecture Notes in Math.}, pages
  233--243. Springer, Berlin, 1979.

\bibitem{kobayashi-vector}
S.~Kobayashi.
\newblock {\em Differential geometry of complex vector bundles}, volume~15 of
  {\em Publications of the Mathematical Society of Japan}.
\newblock Princeton University Press, Princeton, NJ, 1987.
\newblock \protect{Kan\^o} Memorial Lectures, 5.

\bibitem{kobayashi-nagano-filtered-II}
S.~Kobayashi and T.~Nagano.
\newblock On filtered {L}ie algebras and geometric structures. {II}.
\newblock {\em J. Math. Mech.}, 14:513--521, 1965.

\bibitem{luna-inv-75}
D.~Luna.
\newblock Sur les orbites ferm\'ees des groupes alg\'ebriques r\'eductifs.
\newblock {\em Invent. Math.}, 16:1--5, 1972.

\bibitem{mumford-GIT}
D.~Mumford, J.~Fogarty, and F.~Kirwan.
\newblock {\em Geometric invariant theory}, volume~34 of {\em Ergebnisse der
  Mathematik und ihrer Grenzgebiete (2) [Results in Mathematics and Related
  Areas (2)]}.
\newblock Springer-Verlag, Berlin, third edition, 1994.

\bibitem{onishchik-vinberg-seminario}
A.~L. Onishchik and {\`E}.~B. Vinberg.
\newblock {\em Lie groups and algebraic groups}.
\newblock Springer Series in Soviet Mathematics. Springer-Verlag, Berlin, 1990.
\newblock Translated from the Russian and with a preface by D. A. Leites.

\bibitem{ottaviani-rat}
G.~Ottaviani.
\newblock Rational homogeneous varieties.
\newblock Notes from a course held in Cortona, Italy, 1995. \texttt
  {http://www.math.unifi.it/} \texttt{ottavian/}\texttt{public.html}.

\bibitem{ottaviani-spinor}
G.~Ottaviani.
\newblock Spinor bundles on quadrics.
\newblock {\em Trans. Amer. Math. Soc.}, 307(1):301--316, 1988.

\bibitem{ramanan-homogeneous}
S.~Ramanan.
\newblock Holomorphic vector bundles on homogeneous spaces.
\newblock {\em Topology}, 5:159--177, 1966.

\bibitem{umemura-homogeneous}
H.~Umemura.
\newblock On a theorem of {R}amanan.
\newblock {\em Nagoya Math. J.}, 69:131--138, 1978.

\bibitem{wang-xiaowei-balance}
X.~Wang.
\newblock Balance point and stability of vector bundles over a projective
  manifold.
\newblock {\em Math. Res. Lett.}, 9(2-3):393--411, 2002.

\end{thebibliography}
\end{document}